\documentclass{amsproc2000}
\usepackage{latexsym}

\copyrightinfo{2000}{American Mathematical Society}

\theoremstyle{definition}

\theoremstyle{remark}

\usepackage{euscript}
\renewcommand{\mathcal}[1]{\EuScript{#1}}

\let\phi=\varphi

\let\da=\partial

\def\R{\mathbb{R}}
\def\chigh{{\raise1.5pt\hbox{$\chi$}}}

\def\tilq{\widetilde q}

\newcommand{\co}{\colon\thinspace} % FROM G and T, 13-11-99

\let\isom=\cong

\let\sol=\bullet

\def\pback#1{\mathbin{\times_{#1}}}

\def\dminus{\raise2pt\hbox{\vrule height1pt width 2ex}\hskip3pt}

\def\plusH{\ \lower 5pt\hbox{${\buildrel {\textstyle +}
\over {\scriptscriptstyle H}}$}\ }
\def\minusH{\ \lower 5pt\hbox{${\buildrel {\textstyle -}
\over {\scriptscriptstyle H}}$}\ }
\def\timesH{\ \lower 4pt\hbox{${\buildrel {\textstyle .}
\over{\scriptscriptstyle H}}$}\ }

\def\plusV{\ \lower 5pt\hbox{${\buildrel {\textstyle +}
\over {\scriptscriptstyle V}}$}\ }
\def\minusV{\ \lower 5pt\hbox{${\buildrel {\textstyle -}
\over {\scriptscriptstyle V}}$}\ }
\def\timesV{\ \lower 4pt\hbox{${\buildrel {\textstyle .}
\over{\scriptscriptstyle V}}$}\ }

\def\plusK{\ \lower 5pt\hbox{${\buildrel {\textstyle +}
\over {\scriptscriptstyle K}}$}\ }
\def\minusK{\ \lower 5pt\hbox{${\buildrel {\textstyle -}
\over {\scriptscriptstyle K}}$}\ }
\def\timesK{\ \lower 4pt\hbox{${\buildrel {\textstyle .}
\over{\scriptscriptstyle K}}$}\ }

\def\hcompo#1#2#3{{\vcenter{\vbox{\hrule height.#2pt\hbox{\vrule width.#2pt
   height#1pt\kern#3pt\vrule width.#2pt\kern#3pt\vrule width.#2pt}
   \hrule height.#2pt}}}}

\def\vcompo#1#2#3{{\vcenter{\vbox{\hrule height.#2pt
   \hbox{\vrule width.#2pt height#3pt\kern#1pt\vrule width.#2pt}
   \hrule height.#2pt
   \hbox{\vrule width.#2pt height#3pt\kern#1pt\vrule width.#2pt}
   \hrule height.#2pt}}}}

\def\dsq{\mathop{\lower1pt\vbox{\hrule height.4pt \hbox
{\vrule width.4pt height.6em
\kern.6em \vrule width.4pt} \hrule height.4pt}}}

\def\dcomp{\mathop{\dsq\hskip-.88em\raise1pt\hbox{$\scriptstyle\nwarrow$}}}

\def\sq{\vbox{\hrule height.4pt \hbox{\vrule width.4pt height1in
\kern1in \vrule width.4pt} \hrule height.4pt}}

\def\ssq{\vbox{\hrule height.4pt \hbox{\vrule width.4pt height.7in
\kern.7in \vrule width.4pt} \hrule height.4pt}}

\def\tsq{\mathop{\lower1pt\vbox{\hrule height.4pt \hbox
{\vrule width.4pt height.7em
\kern.7em \vrule width.4pt} \hrule height.4pt}}}

\def\sgpd{\,\lower1pt\hbox{$\mlra$}\hskip-0.4in\raise2pt\hbox{$\mlra$}\,}

\def\vgpd{\Bigg\downarrow\!\!\Bigg\downarrow}

\def\surj{-\!\!\!-\!\!\!-\!\!\!\gg}
\def\inj{>\!\!\!-\!\!\!-\!\!\!-\!\!\!>}

\renewcommand{\to}{\longrightarrow}
\let\mapsto=\longmapsto

\def\longmapstoh{\mapstochar\longrightarrow}

\def\longdownmapsto{\raise4ex\hbox{$\Big\downarrow$}}

\def\longdownmaps2{\makebox[0cm]{--}\makebox[0cm]{\raise-2.45ex\hbox{$\Big\downarrow$}}}
\def\gpd{\,\lower1pt\hbox{$\longrightarrow$}\hskip-.24in\raise2pt
             \hbox{$\longrightarrow$}\,}

\def\lrah{\hbox{$\,-\!\!\!-\!\!\!-\!\!\!-\!\!\!-\!\!\!-\!\!\!-\!\!\!\longrightarrow\,$}}
\def\mlra{\hbox{$\,-\!\!\!-\!\!\!\longrightarrow\,$}}

\def\vlra{\hbox{$\,-\!\!\!-\!\!\!-\!\!\!-\!\!\!-\!\!\!-\!\!\!-\!\!\!-\!\!\!
        -\!\!\!-\!\!\!\longrightarrow\,$}}
\def\vgpd{\Bigg\downarrow\!\!\Bigg\downarrow}

\begin{document}

\title[Canonical diffeomorphisms in symplectic 
and Poisson geometry]{On certain canonical diffeomorphisms in symplectic 
and Poisson geometry}

% Remove or comment out any unused author tags.
% author one information

\author{K. C. H. Mackenzie}
\address{Department of Pure Mathematics\\
        University of Sheffield\\
        Sheffield, S3 7RH\\
        United Kingdom}
\email{K.Mackenzie@sheffield.ac.uk}
\thanks{}

% Use this \subjclass if you are using amsproc version 2.0 (December 1999).
\subjclass[2000]{Primary 58H01. Secondary 17B66, 18D05, 22A22, 58F05}
% Use this one if you are using an older version of amsproc.
%\subjclass{}
\date{February 15, 2002}

\begin{abstract}
Associated with the canonical symplectic structure on a cotangent
bundle $T^*M$ is the diffeomorphism $\#\co T^*(T^*M)\to T(T^*M)$. 
This and the Tulczyjew diffeomorphism $T(T^*M)\to T^*(TM)$ may be
derived from the canonical involution $T(TM)\to T(TM)$ by suitable
dualizations. We show that the constructions which yield these maps 
extend very generally to the double Lie algebroids of double Lie 
groupoids, where they play a crucial role in the relations between
double Lie algebroids and Lie bialgebroids. 
\end{abstract}

\maketitle

There have been several talks this meeting about notions of double for Lie
bialgebroids. Some of these have derived from the 1997 construction of 
Liu Zhang--Ju, Alan Weinstein and Xu Ping \cite{LiuWX:1997} in which they 
introduced the notion of Courant algebroid, and some have involved
elements of super mathematics. I very much hope that before
long there will be a clear account of the relations between these various
approaches and even a unification of them. 

There is another approach to the question of doubles, which was not at first
related to Lie bialgebroids, but arose out of broad considerations of what
may be called ``second--order geometry''. It is not possible to describe this 
approach from scratch in an hour, but it is appropriate at this conference 
to indicate the broad features of it. This talk therefore takes a slice
through the papers \cite{Mackenzie:1992}, \cite{Mackenzie:2000}, 
\cite{Mackenzie:1999}, \cite{Mackenzie:notions}, transverse to their 
chronological sequence and provides an alternative route to approach
them. Some aspects of \S2 and \S4 are new. 

I am very grateful to Ted Voronov and Mike Prest for the splendid opportunities 
and good fellowship which the Workshop provided, and to the London 
Mathematical Society for its support. I also wish to thank Yvette
Kosmann--Schwarzbach for her comments on an earlier version. 

\section{}
\label{sect1}

The three canonical diffeomorphisms with which the paper begins are associated 
with iterated tangent and cotangent bundles:

\medskip

{\bf (1)} The canonical involution on an iterated tangent bundle
$$
J\co T^2M\to T^2M,\qquad 
$$
which can be loosely described as `interchanging the order of differentiation'. 

Regard elements of $T^2M = T(TM)$ as represented by equivalence classes
of squares $\phi\co I^2\to M$. Thus $\phi = \phi(s,t)$ differentiates to
a curve of tangent vectors $X_t = \frac{\da\phi}{\da s}(0,t)$ and $X_t$
in turn differentiates to $\xi = \frac{d}{dt}X_t|_0\in T_{X_0}(TM)$. Now
$J$ is induced by the map $\phi\mapsto \phi\circ\sigma$ where 
$\sigma(s,t) = (t,s)$. 

\medskip

{\bf (2)} The Poisson anchor for the canonical symplectic structure on a
cotangent bundle $T^*M$:
$$
\pi^\#_{d\theta}\co T^*(T^*M)\to T(T^*M). 
$$
Any Poisson structure on a manifold $P$ induces a map $\#\co T^*P\to TP$
which sends an exact differential $\delta f$ to the Hamiltonian vector
field $X_f$. This map serves as the anchor of the Lie algebroid structure 
defined by the bracket of 1--forms. When $P$ is the cotangent
bundle $c \co T^*M\to M$, the canonical symplectic structure is 
$\omega = \delta\theta$ where $\theta$ is the 1--form which associates 
to $\eta\in T_\alpha(T^*M)$ the pairing of $\alpha$ with $T(c)(\eta)\in TM$. 
 
\medskip

{\bf (3)} The canonical diffeomorphism 
$$
\Theta\co T(T^*M)\to T^*(TM)
$$
defined intrinsically by Tulczyjew in \cite[p.97]{Tulczyjew} and locally
in, for example, \cite[p.424]{AbrahamM}. We call this the \emph{Tulczyjew
diffeomorphism.} Giving $T(T^*M)$ the tangent lift of the symplectic 
structure on $T^*M$, this $\Theta$ is a symplectic diffeomorphism to the
canonical structure on $T^*(TM)$. 

\bigskip

Anyone who meets these diffeomorphisms for the first time must surely 
suspect that they are related. For example, one would like (2) and (3) 
to be duals of (1). However, it is not at first clear what one might mean 
by this. 

The best way to understand these maps is to extend them into a more general
framework. 

\section{}
\label{sect2}

In this section we consider those aspects of the constructions which do
not depend upon the bracket of vector fields or the dual symplectic
structure. 

First consider $T^2M$. It has two vector bundle structures, the standard 
structure as the tangent bundle of $TM$, for which the projection is
denoted $p_{TM}$, and what we call the \emph{prolongation} structure,
with projection $T(p_M)$ and operations which are obtained by applying the
tangent functor $T$ to the operations in $TM$. These two structures
are represented in the diagram:
$$
\begin{matrix}
	&&T(p_M)&&\\
	&T^2M&\mlra &TM&\cr
        &&&&\cr
 p_{TM} &\Bigg\downarrow&&\Bigg\downarrow&\cr
        &&&&\cr
        &TM&\mlra &M&\cr
\end{matrix}
$$
It is easier to keep the role of the two structures clear in 
one's mind if one replaces $TM$ by an arbitrary vector bundle $(E,q,M)$:
\begin{equation}
\label{TE}
\begin{matrix}
	&&T(q)&&\\
	&TE&\mlra &TM&\cr
        &&&&\cr
 p_{E}  &\Bigg\downarrow&&\Bigg\downarrow&\cr
        &&&&\cr
        &E&\mlra &M&\cr
\end{matrix}
\end{equation}

{\bf Question:} If we dualize $TE\to E$, is there a second structure which
completes the square in
$$
\begin{matrix}
	&T^*E&  & ? &\cr
        &&&&\cr
        &\Bigg\downarrow&&&\cr
        &&&&\cr
        &E&\mlra &M&\cr
\end{matrix}
$$

To see the answer, consider the two short exact sequences associated with $TE$:
\begin{gather*}
q^!E\inj TE \surj q^!TM\\
p_M^! \inj TE \surj p_M^!E
\end{gather*}
Here shrieks denote pullback bundles. The first sequence is a short
exact sequence of vector bundles over base $E$ and the surjection
is the map induced by $T(q)$. The second sequence is over $TM$ and the surjection
is induced by $p_E$. The two surjections are actually the same map 
$TE\to E\times_M TM$, but are considered with respect to different 
bundle structures. 

The kernel of the first surjection is the vertical bundle of $q$; that is, 
the bundle of vectors tangent to the fibres of $E$. Such a vector is determined
by its base point and by a vector in $E$ itself to which it is parallel, and
the injection $q^!E\to TE$ maps $(e_1, e_2)$ to the vertical tangent
vector with base point $e_1$ and parallel to $e_2$. 

The kernel of the second surjection consists of the tangent vectors to $E$
along the zero section. Such a vector $\xi$ may be written as $T(0)(x) + e$
where $x = T(q)(\xi)\in T_m(M)$ and $e\in E_m$ is identified with the vertical 
vector based at $0_m$ and parallel to $e$. 

Thus the $E$ which occurs in both kernels is precisely the set of tangent 
vectors to $E$ which are vertical and based at points of the zero section; that 
is, it is $\ker(T(q))\cap\ker(p_E).$

Dualizing the first short exact sequence gives
$$
q^!T^*M\inj T^*E \surj q^!E^*
$$
and this suggests that the answer to the question is  
\begin{equation}
\label{T^*E}
\begin{matrix}
	&T^*E&\mlra  & E^* &\cr
        &&&&\cr
        &\Bigg\downarrow&&\Bigg\downarrow&\cr
        &&&&\cr
        &E&\mlra &M&\cr
\end{matrix}
\end{equation}
That this is indeed so was proved in \cite{MackenzieX:1994}. 

To find the dual of the other structure, $TE\to TM$, apply $T$ to the
pairing $E\times_M E^*\to\R$. This gives a pairing
$TE\times_{TM}T(E^*)\to\R$ which is again non--degenerate. So
\begin{equation}
\label{sol}
(TE\to TM)^*\isom (T(E^*)\to TM)
\end{equation}
and so dualizing the horizontal structure in (\ref{TE}) gives
\begin{equation}
\label{TE^*}
\begin{matrix}
	&T(E^*)&\mlra  & TM &\cr
        &&&&\cr
        &\Bigg\downarrow&&\Bigg\downarrow&\cr
        &&&&\cr
        &E^*&\mlra &M&\cr
\end{matrix}
\end{equation}
This is precisely the structure (1) with $E$ replaced by $E^*$. 

Now consider more general structures of this type. 

A \emph{double vector bundle} (\cite{Pradines:DVB}; see also 
\cite[\S1]{Mackenzie:1992}) is a manifold $D$ with two vector bundle
structures, on bases $A$ and $B$, each of which is a vector bundle on
base $M$,
$$
\begin{matrix}
	&&\tilq_H&&\\
	&D&\mlra &B&\cr
        &&&&\cr
\tilq_V &\Bigg\downarrow&&\Bigg\downarrow&q_B\cr
        &&&&\cr
        &A&\mlra &M&\cr
	&&q_A&&\
\end{matrix}
$$
such that the structure maps of each vector bundle structure on $D$ are
vector bundle morphisms with respect to the other. Precisely: for the 
horizontal structure,

$\sol$ the projection $\tilq_H\co D\to B$ is a vector bundle morphism over the
projection $q_A\co A\to M$;

$\sol$ the zero section $B\to D$ is a vector bundle morphism over the
zero section $M\to A$;

$\sol$ the addition $D\pback{B}D\to D$ is a vector bundle morphism
over the addition $A\pback{M}A\to A$;

$\sol$ the scalar multiplication $D\times\R\to D$ is a vector bundle morphism
over the scalar multiplication $A\times\R\to A$.

\medskip

The first three conditions are precisely those which ensure that for four 
elements $d_1, d_2, d_3, d_4\in D$ such that the LHS of
$$
(d_1 \plusH d_2)\plusV	(d_3 \plusH d_4) =
(d_1\plusV d_2)\plusH (d_3\plusV d_4)
$$
is defined, the RHS is defined and the two sides are equal.
That (\ref{TE}) is indeed a double vector bundle is easily verified. 

Now consider the intersection of the two kernels
$$
C = \ker(\tilq_H)\cap\ker(\tilq_V).
$$
A priori, there is no reason to expect this to have a natural algebraic
structure. However, each structure on $D$ induces a vector bundle structure
on $C$ with base $M$, and these two structures coincide. Pradines'
terminology in \cite{Pradines:DVB} was \emph{c\oe ur}; we call it the 
\emph{core} of $D$ . 

There is again a short exact sequence
$$
q_A^!C\inj D\surj q_A^!B
$$
of vector bundles over $A$ and dualizing this gives
$$
q_A^!B^*\inj D^{*V}\surj q_A^!C^*
$$
and we get another double vector bundle
$$
\begin{matrix}
	&D^{*V}&\mlra &C^*&\cr
        &&&&\cr
        &\Bigg\downarrow&&\Bigg\downarrow&\cr
        &&&&\cr
        &A&\mlra &M&\cr
\end{matrix}
$$
called the \emph{vertical dual} of $D$. Likewise there is a \emph{horizontal
dual}
$$
\begin{matrix}
	&D^{*H}&\mlra &B&\cr
        &&&&\cr
        &\Bigg\downarrow&&\Bigg\downarrow&\cr
        &&&&\cr
        &C^*&\mlra &M&\cr
\end{matrix}
$$

{\bf Theorem 1:} (\cite{Mackenzie:1999}, \cite{KoniecznaU:1999})
The vector bundles $D^{*V}\to C^*$ and $D^{*H}\to C^*$ are themselves
dual under a pairing which is canonical up to sign and which respects
the double vector bundle structures. 

\medskip

To define the pairing, consider elements $\Phi\in D^{*V}$ and $\Psi\in D^{*H}$, 
shown with their various projections as follows:
$$
\begin{matrix}
          \Phi&\longmapstoh&\kappa\cr
          \longdownmaps2&&\longdownmaps2\cr
          &&\cr
          a&\longmapstoh&m\cr
\end{matrix}
          \qquad\qquad\qquad
\begin{matrix}
          \Psi&\longmapstoh&b\cr
          \longdownmaps2&&\longdownmaps2\cr
          &&\cr
          \kappa&\longmapstoh&m\cr
\end{matrix}
$$
where $a\in A,\ b\in B$ and $\kappa\in C^*$. The given elements must of course
lie over the same point of $C^*$. Now given $a$ and $b$
over the same $m\in M$, there exist elements $d\in D$ of the form
$$
\begin{matrix}
          d&\longmapstoh&b\cr
          \longdownmaps2&&\longdownmaps2\cr
          &&\cr
          a&\longmapstoh&m\cr
\end{matrix}
$$
and we define
$$
\langle\Phi, \Psi\rangle_{C^*} = 
\langle\Phi, d\rangle_{A} - \langle\Psi, d\rangle_{B} 
$$
where the subscripts indicate the bases of the pairings. 

We could equally well reverse the order on the RHS. The statement in the theorem
that this duality respects the double vector bundle structures means that
the induced isomorphism from the dual of $D^{*V}\to C^*$ to $D^{*H}\to C^*$
is in fact an isomorphism of double vector bundles. 

\medskip

{\bf Example:} 
Consider $D = TE$ as in (\ref{TE}). We observed that the two duals are
$D^{*V} = T^*E$ as in (\ref{T^*E}) and $D^{*H} = T(E^*)$ as in
(\ref{TE^*}). So Theorem 1 implies that there is a non--degenerate pairing
of the two bundles
$$
T^*E\to E^* \qquad \mbox{and}\qquad T(E^*)\to E^*. 
$$
To recognize what this is, write $A = E^*$. Then the statement is that the two
bundles
$$
T^*(A^*)\to A\qquad \mbox{and}\qquad T^*(A)\to A
$$
are isomorphic as double vector bundles. That this is so was proved (by
different methods) in \cite{MackenzieX:1994}. To state the result precisely: 
Given a vector bundle $A\to M$, there is a canonical isomorphism of double 
vector bundles $R_A\co T^*(A^*)\to T^*(A)$ over $A$ and $A^*$, which is an 
anti--symplectomorphism and which induces the negative of the identity on 
the cores $T^*M\to T^*M$. 

A proof of the next result is also in \cite{MackenzieX:1994}. 

\medskip

{\bf Proposition:}
For any manifold $M$, the composite
$$
T^*(T^*M) \stackrel{\pi^\#_{d\theta}}{\longrightarrow} T(T^*M)
\stackrel{\Theta}{\longrightarrow} T^*(TM)
$$
is equal to $R_{TM}$. 

\medskip

The significance of this is that $R$ is defined for arbitrary vector bundles, 
without reference to any structures of bracket type. Thus the Tulczyjew 
isomorphism could be intrinsically defined in terms of the symplectic structure
on $T^*M$ and the map $R$. (The definition in \cite{Tulczyjew} uses the
canonical involution on $T^2M$ and the pairing which induces (\ref{sol}); 
see the example at the end of \S\ref{sect3} below.)

Before we move on to consider the role of bracket structures in this context, 
it is interesting to observe that Theorem 1 is forced upon us by the
simple pursuit of diagrams. Take an arbitrary double vector bundle and 
apply the tangent functor to get the triple structure in Figure~\ref{fig:1}(a). 
\begin{figure}[h] % [b]
\begin{picture}(340,180)(0,30)               % \label{pic:triples}
\put(-10,150){$\begin{matrix}
               &&     &&\cr
                     &TD&\mlra & TB & \cr
                     &&&&\cr
                     &\Bigg\downarrow  & &\Bigg\downarrow  &   \cr
                     &&&&\cr
                     &D&\mlra &B&\cr
               \end{matrix}$}

\put(30, 165){\vector(1,-1){20}}                % Top left
\put(100, 165){\vector(1,-1){18}}               % Top right
\put(30, 100){\vector(1,-1){20}}                % Bottom left
\put(100, 100){\vector(1,-1){20}}               % Bottom right

\put(38,105){$\begin{matrix}
               &&      &\cr
                     &TA &\mlra &TM\cr
                     &&&\cr
                     &\Bigg\downarrow &&\Bigg\downarrow \cr
                     &&&\cr
                     &A &\mlra & M\cr
               \end{matrix}$}

\put(90,20){(a)}

% RIGHT HAND HALF

\put(210,150){$\begin{matrix}
                &&      &\cr
                      &T^*D&\mlra &D^{*H}\cr
                      &&&\cr
                      &\Bigg\downarrow & &\Bigg\downarrow \cr
                      &&&\cr
                      &D&\mlra &B \cr
                \end{matrix}$}

\put(245, 165){\vector(1,-1){20}}                % Top left
\put(320, 165){\vector(1,-1){20}}               % Top right
\put(245, 100){\vector(1,-1){20}}                % Bottom left
\put(320, 100){\vector(1,-1){20}}               % Bottom right

\put(260,105){$\begin{matrix}
               &&      &\cr
                     &D^{*V}&\mlra &C^*\cr
                     &&&\cr
                     &\Bigg\downarrow &&\Bigg\downarrow \cr
                     &&&\cr
                     &A &\mlra & M\cr
               \end{matrix}$}
\put(310,20){(b)}
\end{picture}\caption{\ \label{fig:1}}
\end{figure}
(The diagrams should be viewed as coming out of the page to the right.)
Now dualize the bundle $TD\to D$. The rear square is precisely of the
form (\ref{TE}), but when we dualize $D$ we do so over $B$ and so the
dual is the horizontal dual. Likewise the left side face is of the form
(\ref{TE}) but when we dualize $D$ we do so over $A$, and so obtain the
vertical dual. Now it is easy to guess that the eighth vertex must be $C^*$
and one can then verify that Figure~\ref{fig:1}(b)
is indeed a triple vector bundle. Using the isomorphisms $R$ for the left
and rear sides, we know that $T^*D\isom T^*(D^{*V})$ and 
$T^*D\isom T^*(D^{*H})$ as double vector bundles and we are forced to the 
conclusion that the top face must also be of the form (\ref{T^*E}), and 
that $D^{*V}$ and $D^{*H}$ must be duals. 

\section{}
\label{sect3}

Bracket structures arise naturally from the differentiation of Lie
groupoids. To proceed further, we thus need to consider general double 
Lie groupoids as in \cite{Mackenzie:1992}. 
\begin{equation}
\label{S}
\begin{matrix}
%        &&  &&\cr
        &S &\sgpd   &V&\cr
        &&&&\cr
        &\vgpd & &\vgpd & \cr
        &&&&\cr
        &H&\sgpd &M&\cr
%        &&&&\cr
%        &&&&\cr
%        &&\mbox{(a)}&&\cr
\end{matrix}
\end{equation}

Applying the Lie functor to the vertical structures gives the diagram on the 
left below \cite{Mackenzie:1992}. That the Lie algebroid $A_VS$ of the 
vertical structure inherits a groupoid structure with base $AV$ follows directly
by observing that the Lie functor --- like the tangent functor --- preserves 
diagrams and pullbacks. The horizontal groupoid structures and the vertical Lie
algebroid structures are compatible in the same way as in the definition
of a double vector bundle above. Now applying the Lie functor a second time, 
to the horizontal structures, gives the diagram on the right 
\cite{Mackenzie:2000}. This is a double Lie algebroid, but for the moment 
we will be concerned with it as a double vector bundle. 
\begin{equation}
\label{A_V}
\begin{matrix}
%        &&  &&\cr
        &A_VS &\sgpd   &AV&\cr
        &&&&\cr
        &\Big\downarrow& &\Big\downarrow& \cr
        &&&&\cr
        &H&\sgpd &M&\cr
%        &&&&\cr
%        &&&&\cr
%        &&\mbox{(a)}&&\cr
\end{matrix}
\hskip2cm
\begin{matrix}
%        &&  &&\cr
        &A_H(A_VS) &\longrightarrow   &AV&\cr
        &&&&\cr
        &\Big\downarrow& &\Big\downarrow& \cr
        &&&&\cr
        &AH&\longrightarrow &M&\cr
%        &&&&\cr
%        &&&&\cr
%        &&\mbox{(a)}&&\cr
\end{matrix}
\end{equation}
We could equally well apply the Lie functor first horizontally and then
vertically, obtaining the following structures. 
\begin{equation}
\label{A_H}
\begin{matrix}
%        &&  &&\cr
        &A_HS &\longrightarrow   &V&\cr
        &&&&\cr
        &\vgpd & &\vgpd & \cr
        &&&&\cr
        &AH&\longrightarrow &M&\cr
%        &&&&\cr
%        &&&&\cr
%        &&\mbox{(a)}&&\cr
\end{matrix}
\hskip2cm
\begin{matrix}
%        &&  &&\cr
        &A_V(A_HS) &\longrightarrow   &AV&\cr
        &&&&\cr
        &\Big\downarrow& &\Big\downarrow& \cr
        &&&&\cr
        &AH&\longrightarrow &M&\cr
%        &&&&\cr
%        &&&&\cr
%        &&\mbox{(a)}&&\cr
\end{matrix}
\end{equation}

The concept of core extends to double groupoids \cite{BrownM:1992},
\cite{Mackenzie:1992} and in the case of $S$ above yields a Lie groupoid 
$K$ on base $M$, the Lie algebroid of which is the core of the two double
vector bundles $A_V(A_HS)$ and $A_H(A_VS)$. 

The Lie functor is a kind of generalized tangent functor, and both 
$A_V(A_HS)$ and $A_H(A_VS)$ are submanifolds of $T^2S$. Further, the 
canonical involution $J$ carries $A_V(A_HS)$ onto $A_H(A_VS)$ 
\cite{Mackenzie:2000}. 
Denote this restricted map by $j_S$. This is an isomorphism of double vector
bundles, as in Figure~\ref{fig:j}.  

\begin{figure}[h]
\begin{picture}(200,180)(60,0)
% \put(0,0){.}                          % checker
% \put(155,120){.}                         % checker
\put(0,120){$\begin{matrix}
                      &&   &\cr
                      &A_H(A_VS)&\lrah &AV\cr
                      &&&\cr
                      &\Bigg\downarrow& &\Bigg\downarrow\cr
                      &&&\cr
                      &AH&\lrah &M\cr
            \end{matrix}
                      $}
\put(60,140){\vector(2,-1){120}}                               % top left
% \put(5,7){\vector(4,-3){3}}                              % top right
% \put(1.1,5){\line(4,-3){3}}                                 % bottom left
% \put(1,5){\line(4,-3){3}}                                 % bottom left

\put(90,110){$j_S$}
% \put(8,5){$\phi^V$}

\put(185,45){$\begin{matrix}
                     &&    &\cr
                     &A_V(A_HS)&\lrah &AV\cr
                     &&&\cr
                     &\Bigg\downarrow&&\Bigg\downarrow\cr
                     &&&\cr
                     &AH&\lrah & M\cr
               \end{matrix}
                     $}
\end{picture}\caption{\ \label{fig:j}}
\end{figure}

Dualizing $j_S$ over $AH$ we obtain a map
$$
j^{*V}\co A_V^*(A_HS)\to A_H^\sol(A_VS). 
$$
The structure on the domain is the ordinary dual of the Lie algebroid of
a Lie groupoid. On the target, however, the solidus indicates that we have
taken the dual of a prolongation structure. As with the case
$$
% T^\sol E \isom 
\begin{pmatrix}
TE\\
\Bigg\downarrow\\
TM
\end{pmatrix}^*
\isom
\begin{pmatrix}
T(E^*)\\
\Bigg\downarrow\\
TM
\end{pmatrix}
$$
in (\ref{sol}), we can show that
$$
% T^\sol E \isom 
\begin{pmatrix}
A_H(A_VS)\\
\Bigg\downarrow\\
AH
\end{pmatrix}^*
\isom
\begin{pmatrix}
A(A_V^*S)\\
\Bigg\downarrow\\
AH
\end{pmatrix}
$$
where 
\begin{equation}
\label{dualvb}
\begin{matrix}
%        &&  &&\cr
        &A_V^*S &\sgpd   &A^*K&\cr
        &&&&\cr
        &\Big\downarrow& &\Big\downarrow& \cr
        &&&&\cr
        &H&\sgpd &M&\cr
%        &&&&\cr
%        &&&&\cr
%        &&\mbox{(a)}&&\cr
\end{matrix}
\end{equation}
is a structure obtained by applying to $A_VS$ an extended form of the duality 
for double vector bundles (see \cite{Pradines:1986} or \cite{Mackenzie:1999}). 

Denote the isomorphism $A(A_V^*S)\to A_H^\sol(A_VS)$ by $I_V$. We then define
a generalization of the Tulczyjew isomorphism by
$$
\Theta = I_V^{-1}\circ j^{*V}\co A_V^*(A_HS)\to A(A_V^*S)
$$
(In \cite{Mackenzie:1999} the inverse of this map is denoted ${j'}^{V}$.) 

\medskip

{\bf Example:} 
To recover the classical case, consider the double groupoid
$$
\begin{matrix}
%        &&  &&\cr
    S = &M^4 &\sgpd   &M^2&\cr
%        &&&&\cr
        &\vgpd & &\vgpd & \cr
%        &&&&\cr
        &M^2&\sgpd &M&\cr
%        &&&&\cr
%        &&&&\cr
%        &&\mbox{(a)}&&\cr
\end{matrix}
$$
in which elements are quadruples of points, regarded as the corners of empty
squares. Applying the Lie functor vertically, we obtain
$$
\begin{matrix}
%        &&  &&\cr
        A_VS = & TM\times TM &\sgpd   &TM&\cr
%        &&&&\cr
        &\Big\downarrow& &\Big\downarrow& \cr
%        &&&&\cr
        &M\times M&\sgpd &M&\cr
%        &&&&\cr
%        &&&&\cr
%        &&\mbox{(a)}&&\cr
\end{matrix}
$$
Here $TM\times TM$ is both the tangent Lie algebroid of $M\times M$ and the
pair groupoid on base $TM$. The dualization process then produces
$$
 \begin{matrix}
%        &&  &&\cr
        A_V^*S = & T^*M\times T^*M &\sgpd   &T^*M&\cr
%        &&&&\cr
        &\Big\downarrow & &\Big\downarrow & \cr
%        &&&&\cr
        &M\times M &\sgpd &M&\cr
%        &&&&\cr
%        &&&&\cr
%        &&\mbox{(a)}&&\cr
\end{matrix}
$$
in which again the horizontal groupoid structure is a pair groupoid. Applying
the Lie functor to this horizontally we obtain
$$
\begin{matrix}
%        &&  &&\cr
        A(A_V^*S) = & T(T^*M) &\longrightarrow   &T^*M&\cr
%        &&&&\cr
        &\Big\downarrow& &\Big\downarrow& \cr
%        &&&&\cr
        &TM&\longrightarrow &M&\cr
%        &&&&\cr
%        &&&&\cr
%        &&\mbox{(a)}&&\cr
\end{matrix}
$$
The reader can similarly check that $A^*_V(A_HS) = T^*(TM)$. In this case
$I_V$ is precisely the isomorphism induced by the tangent pairing (\ref{sol}) 
of $T^2M$ and $T(T^*M)$, and the definition of $\Theta$ in this
case thus reduces exactly to that in \cite{Tulczyjew}. 

\section{}
\label{sect4}

Finally we consider a general form of the map $T^*(T^*M)\to T(T^*M)$
obtained from the symplectic structure on a cotangent bundle. 

\medskip

{\bf Theorem 2:} \cite{Mackenzie:1999} Given any double Lie groupoid (\ref{S})
with core groupoid $K\gpd M$, there is a
structure of double Lie groupoid on $T^*S$ which makes it
a symplectic double groupoid
$$
\begin{matrix}
%        &&  &&\cr
        &T^*S &\sgpd   &A_V^*S&\cr
        &&&&\cr
        &\vgpd & &\vgpd & \cr
        &&&&\cr
        &A_H^*S&\sgpd &A^*K&\cr
%        &&&&\cr
%        &&&&\cr
%        &&\mbox{(a)}&&\cr
\end{matrix}
$$
and which induces on $A_V^*S$, $A_H^*S$ and $A^*K$ their standard Poisson 
structures. In particular, 
$A_H^*S\gpd A^*K$ and $A_V^*S\gpd A^*K$ are Poisson groupoids in duality. 

\medskip

The statement that 
% $A_H^*S\gpd A^*K$ and $A_V^*S\gpd A^*K$ are 
two Poisson groupoids are in duality means that the 
Lie algebroid of each, constructed from the Lie groupoid structure, is 
isomorphic (or anti--isomorphic) to the Lie algebroid dual of the other 
with the Lie algebroid structure which it acquires from the Poisson 
groupoid structure \cite{Weinstein:1988}. (We will not consider how to arrange 
signs here; there are various possibilities.) Thus Theorem 2 asserts the 
existence of Lie algebroid isomorphisms
$$
\mathcal{D}_H\co A^*(A_H^*S)\to A(A_V^*S),\qquad
\mathcal{D}_V\co A^*(A_V^*S)\to A(A_H^*S). 
$$
Following through the Example at the end of \S\ref{sect3}, the reader can verify 
that these generalize the map $\pi^\#_{d\theta}$. Whereas $\pi^\#_{d\theta}$
arises from the symplectic structure on $T^*M$, the maps $\mathcal{D}_H$ and 
$\mathcal{D}_V$ arise from the Poisson structures on the duals of Lie algebroids. 

In order to express the relationship between these maps $\mathcal{D}$ and
the Tulczyjew diffeomorphism in the general case, we need one more ingredient. 
Consider the structure $A_HS$ from (\ref{A_H}). Applying the duality process
mentioned in connection with (\ref{dualvb}), we obtain the structure on the
left below. Applying the Lie functor vertically we obtain the double vector
bundle on the right below. 
\begin{equation}
\label{hmmm}
\begin{matrix}
%        &&  &&\cr
        &A_H^*S &\longrightarrow   &V&\cr
        &&&&\cr
        &\vgpd & &\vgpd & \cr
        &&&&\cr
        &A^*K&\longrightarrow &M&\cr
%        &&&&\cr
%        &&&&\cr
%        &&\mbox{(a)}&&\cr
\end{matrix}
\hskip2cm
\begin{matrix}
%        &&  &&\cr
        &A(A_H^*S) &\longrightarrow   &AV&\cr
        &&&&\cr
        &\Big\downarrow& &\Big\downarrow& \cr
        &&&&\cr
        &A^*K&\longrightarrow &M&\cr
%        &&&&\cr
%        &&&&\cr
%        &&\mbox{(a)}&&\cr
\end{matrix}
\end{equation}
Now consider the vertical dual of the double vector bundle $A_V(A_HS)$ in
(\ref{A_H}). By an extension of the methods used in \S\ref{sect2}, one can
define a non--degenerate pairing of the double vector bundles $A(A_H^*S)$ 
and $A_V^*(A_HS)$which induces an isomorphism of double vector bundles
$$
R_H\co A_V^*(A_H^*S)\to A^*(A_HS).
$$
In the case where $S$ is actually a double vector bundle, this reduces to
the $R$ of \S\ref{sect2}. 

\medskip

{\bf Theorem 3:}\cite{Mackenzie:1999}
For a double Lie groupoid $S$ as in (\ref{S}), the following diagram of
diffeomorphisms commutes. 

\begin{picture}(200,150)
% \put(0,0){.}                          % checker
% \put(155,120){.}                         % checker
\put(0,60){$\begin{matrix}
       &&\mathcal{D}_H^{-1}  &&\cr
       &A(A_V^*S) &\vlra   &A^*(A_H^*S)&\cr
       &&&&\cr
I_V    &\Big\downarrow& &\Big\downarrow&R_H \cr
       &&&&\cr
       &A^\sol(A_HS)&\vlra   &A^*(A_HS)&\cr
       &&(j^{*V})^{-1}&&\cr
%        &&&&\cr
%        &&\mbox{(a)}&&\cr
\end{matrix}$}
\put(70,80){\vector(2,-1){80}}                      
\put(100,70){$\Theta$}
\end{picture}

\medskip

In particular, it follows that $\Theta$ is a Poisson diffeomorphism, since 
both $\mathcal{D}_H$ and $R_H$ are anti--Poisson \cite{Mackenzie:1999}. 

\bigskip

In conclusion, it is worthwhile noticing the different degrees of 
development required by the various canonical diffeomorphims which we have 
discussed. The map $R_A\co T^*(A^*)\to T^*(A)$ is the simplest, requiring 
only a single ordinary vector bundle. The Tulczyjew diffeomorphism requires 
the notion of double Lie groupoid, and of its associated double Lie algebroid, 
but involves only the simpler aspects of that construction --- a dualization 
of the map $j$ as a morphism of double vector bundles and the prolongation 
of a pairing which exists on the level of the intermediate structures $A_VS$. 
The maps $\mathcal{D}_V$ and $\mathcal{D}_H$, on the other hand, require a 
detailed consideration of the structure of the symplectic double groupoid 
associated with a double Lie groupoid. As was shown in \cite{Mackenzie:1999}
and \cite{Mackenzie:notions}, this structure is the key to the general
notion of double Lie algebroid. 

% following is corect for pms2 and main
% \bibliography{/home/main/BIB/pream,/home/main/BIB/lgla,%
% /home/main/BIB/general,/home/main/BIB/preprint}
% \bibliographystyle{plain}

\newcommand{\noopsort}[1]{} \newcommand{\singleletter}[1]{#1}

\end{document}